\documentclass[onecolumn]{autart}    

\usepackage{graphicx}      
\usepackage{amsmath,amsfonts, amssymb,subcaption, cite}
\usepackage{tikz}
\usetikzlibrary{shapes.geometric}

\newtheorem{corollary}{Corollary}
\newtheorem{proposition}{Proposition}
\newtheorem{definition}{Definition}

\DeclareMathOperator{\diag}{diag}

\begin{document}

\begin{frontmatter}

\title{A Graph-Theoretic Approach to the $\mathcal{H}_{\infty}$ Performance of Dynamical Systems on  Directed and Undirected Networks\thanksref{footnoteinfo}}

\thanks[footnoteinfo]{The material in this paper was not presented at any conference.}

\author[First]{Mohammad Pirani},
\author[First]{Henrik Sandberg}, 
\author[First]{Karl Henrik Johansson}

\address[First]{Department of Automatic Control, KTH Royal Institute of Technology, Stockholm, Sweden (e-mail: pirani, hsan,
kallej@kth.se)}


\begin{abstract}                           We study a graph-theoretic approach to the $\mathcal{H}_{\infty}$ performance of leader following consensus dynamics on directed and undirected graphs. We first provide graph-theoretic bounds on the system $\mathcal{H}_{\infty}$ norm of the leader following dynamics and show the tightness of the proposed bounds. Then, we discuss the relation between the system $\mathcal{H}_{\infty}$ norm for directed and undirected networks for specific classes of graphs, i.e., balanced digraphs and directed trees. Moreover, we investigate the effects of adding directed edges to a directed tree on the resulting system $\mathcal{H}_{\infty}$ norm. In the end, we apply these theoretical results to  a reference velocity tracking problem in a platoon of connected vehicles and discuss the effect of the location of the leading vehicle on the overall $\mathcal{H}_{\infty}$ performance of the system.
\end{abstract}

\end{frontmatter}

\section{Introduction}\label{sec:intro} 

The notion of robustness of dynamical systems to external disturbances or parameter uncertainties has been under investigations for the past decades \cite{Doyle}. The theory of robust control started to shape in the 1970s and 1980s when the lack of robustness in the state-space representation of dynamical systems was felt. 
The developments of optimization techniques enabled robust control problems to be recast into semidefinite programming or eigenvalue problems which could  be  efficiently solved with computers \cite{Boyyd, Dullerud}. Among optimization based methods, $\mathcal{H}_{\infty}$ techniques became popular in handling the design of robust MIMO systems \cite{Dolyfrancis}.

By enlarging the size of communication networks for various applications and
introducing new interdisciplinary areas such as the Internet of Things and Cyber-Physical Systems, the need to address the safety and robustness of large-scale networks was seriously felt. In this direction, as the size of the plant (network) increases and the interactions become more sophisticated, having a knowledge about the  the network structure can help us to address the resilience and robustness of such complex systems more efficiently. In this direction, there has been a vast literature in analyzing the effect of network structure on both system $\mathcal{H}_{2}$, \cite{Chapman, Siamii, Scardovi} , and $\mathcal{H}_{\infty}$ performances \cite{ArxiveRobutness, Wei, jay2}. Via combining the system-theoretic notions with algebraic graph theory, some papers have looked at these performance metrics as network centrality measures and discussed the control node (leader) selection problems in a given large-scale network to optimize each performance metric \cite{Fitch2, ArxiveRobutness,  Siami2}.

Since for a linear-time-invariant system, the $\mathcal{H}_{\infty}$ norm is the worst-case gain over all frequencies, it may exhibit specific and contradictory behaviors compared to other system norms \cite{PiraniJohnCDC}. On the other hand, when the interactions between agents in a communication network are not symmetric, various (and probably counter-intuitive) behaviors can be expected from the dynamics defined on this network \cite{Scardovi2}. The current paper analyzes the $\mathcal{H}_{\infty}$  performance of dynamical systems on directed networks. More specifically, the contributions of this paper are: 
\begin{enumerate}
    \item We discuss some graph-theoretic bounds on the system $\mathcal{H}_{\infty}$ norm on directed and undirected graphs and show the tightness of the proposed bounds via examples. 
    
    \item We investigate the relation between system $\mathcal{H}_{\infty}$  norms in directed and undirected networks for specific classes of networks, i.e., balanced digraphs and directed trees. Moreover, we discuss the effect of adding (or removing) directed edges on the $\mathcal{H}_{\infty}$  performance of the system. 
    
    \item We apply these results to discuss the effect of the directed network and the location of the leading vehicle on a reference velocity tracking scenario in vehicle platooning.

\end{enumerate}

After introducing some mathematical notations and definitions in Section \ref{sec:not}, in Section \ref{sec:influenawce} we state the problem of the $\mathcal{H}_{\infty}$ performance of a reference-following network dynamical system. In Section \ref{sec:boundssss}, graph-theoretic  bounds on the system $\mathcal{H}_{\infty}$ norms for undirected and directed networks are introduced and based on that, necessary and sufficient conditions for a network system to have $\mathcal{H}_{\infty}$ norm less than a specific number is discussed. In Section \ref{sec:relations} we compare the $\mathcal{H}_{\infty}$ norm of a directed network with its undirected counterpart for balanced digraphs and directed trees and the effect of adding directed edges to a directed tree on the $\mathcal{H}_{\infty}$ performance of the dynamics. In Section \ref{sec:platoonn} we apply some of the results to discuss the robustness of a vehicle platooning scenario. Section \ref{sec:conclusion} concludes the paper.

\section{Notations and Definitions}
\label{sec:not}

We use $\mathcal{G}_d=\{\mathcal{V},\mathcal{E}\}$ to denote an unweighted directed graph where $\mathcal{V}$ is the set of vertices (or nodes) and $\mathcal{E}$ is the set of directed edges, i.e., $(v_i,v_j)\in \mathcal{E}$ if an only if there exists a directed edge from $v_i$ to $v_j$. Moreover, an unweighted undirected graph $\mathcal{G}_u=\{\mathcal{V},\mathcal{E}\}$ is a graph such that  $(v_i,v_j)\in \mathcal{E}$ if an only if there exists an undirected edge between $v_i$ and $v_j$.  For directed graphs in this paper, we only consider unidirectional edges, i.e., if there exists a direct edge $v_i$ to $v_j$, then there is no direct edge from $v_j$ to $v_i$. Let $|\mathcal{V}|=n$ and define the adjacency matrix for $\mathcal{G}_d$, denoted by $A_{n\times n}$, to be a binary matrix  where  $A_{ij}=1$ if and only if there is a directed edge from $v_j$ to $v_i$ in $\mathcal{G}_d$. The {\it neighbors} of vertex $v_i \in \mathcal{V}$ in graph $\mathcal{G}_d$ are given by the set $\aleph_i = \{v_j \in \mathcal{V}~|~(v_j, v_i) \in \mathcal{E}\}$. We define the in-degree for node $v_i$ as $\Delta_i=\sum_{v_j\in \aleph_i} A_{ij}$ and the out-degree as $\delta_i=\sum_{v_j\in \aleph_i} A_{ji}$. For a given set of nodes $X \subset \mathcal{V}$, the {\it edge-boundary} (or just boundary) of the set is defined as $\partial{X}= \{(v_i,v_j) \in \mathcal{E} | v_i \in X, v_j \in \mathcal{V}\setminus{X}\}$. For a symmetric matrix $M$, the eigenvalues are ordered as  $\lambda_1(M) \leq \lambda_2(M) \leq \ldots \leq \lambda_n(M)$ and the singular values of a matrix $\mathcal{M}$ are ordered as $\sigma_1(\mathcal{M}) \leq \sigma_2(\mathcal{M}) \leq ... \leq \sigma_n(\mathcal{M})$. The Laplacian matrix of the graph is $\mathcal{L} = D - A$, where $D = \diag (\Delta_1, \Delta_2, ..., \Delta_n)$. We will be considering a nonempty subset of  vertices $\mathcal{\mathfrak{S}} \subset \mathcal{V}$ to be {\it leaders}, whose in-degree is zero, and assume without loss of generality that the leaders are placed last in an ordering of the agents. Vertices in $\mathcal{V}\setminus\mathfrak{S}$ are called {\it followers}. The {\it grounded Laplacian} induced by the leader set $\mathfrak{S}$ is denoted by $\mathcal{L}_g(\mathfrak{S})$ or simply $\mathcal{L}_g$, and is obtained by removing the rows and columns of $\mathcal{L}$ corresponding to the nodes in $\mathfrak{S}$ \cite{ACC}. If the underlying graph is directed, the grounded Laplacian is denoted by $\mathcal{L}_{g,d}$ and if the graph is undirected, it is $\mathcal{L}_{g,u}$. 
The state-space representation of a linear time-invariant system with $n$ states, $m$ inputs and $k$ outputs is denoted by the triple $(A_{n\times n}, B_{n\times m}, C_{k \times n})$, where $A$ is the state matrix, $B$ is the input matrix and $C$ is the output matrix.  We use $\mathbf{e}_i$ to indicate the $i$-th vector of the canonical basis.

\section{Problem Statement}
\label{sec:influenawce}

Consider a connected network  consisting of $n$ agents $\mathcal{V} = \{v_1, v_2,\ldots, v_n\}$. The set of agents is partitioned into a set of followers $\mathfrak{F}$, and a set of leaders\footnote{These agents may also be referred to as {\it{anchors}} \cite{Rahmani} or {\it{stubborn agents}}\cite{Ghaderi13} depending on the context.} $\mathfrak{S}$. We assume that there exists at least one leader in the network. The number of leaders which are connected to follower $v_i$ is denoted by $\Gamma_i$, i.e., $\Gamma_i\triangleq |\aleph_i\cap \mathfrak{S}|$,  and we denote $\Gamma_{\rm max}=\max_{i\in \mathcal{V}\setminus \mathfrak{S}}\Gamma_i$ and $\Gamma_{\rm min}=\min_{i\in \mathcal{V}\setminus \mathfrak{S}}\Gamma_i$. Examples of such directed graphs are shown in Fig.~\ref{fig:etfdvnj}, where black nodes are leaders and white nodes are the followers. Each agent $v_i$ has a scalar and real valued state $\psi_i(t)$, where $t$ is the time index. The state of each follower agent $v_j\in \mathfrak{F}$ evolves based on the interactions  with its neighbors as
\begin{align}
\dot{\psi}_{j}(t)&=\sum_{v_i\in \aleph_j}(\psi_i(t)-\psi_j(t)).
\label{eqn:partial}
\end{align}
Since the leaders state are not influence by the followers, their state  is assumed to be constant and thus
\begin{equation}
\dot{\psi}_{j}(t) = 0, \enspace \forall v_j \in \mathfrak{S}.
\label{eqn:fully}
\end{equation}
 Aggregating the states of all followers  into a vector $\boldsymbol \psi_\mathfrak{F}(t) \in \mathbb{R}^{n-|\mathfrak{S}|}$, and the states of all  leaders into a vector $\boldsymbol \psi_{\mathfrak{S}}(t)\in \mathbb{R}^{|\mathfrak{S}|}$ (note that $\boldsymbol \psi_{\mathfrak{S}}(t) = \boldsymbol \psi_{\mathfrak{S}}(0)$ for all $t \ge 0$), equations \eqref{eqn:partial} and \eqref{eqn:fully} yield the following dynamics 
\begin{equation}
\begin{bmatrix}
      \dot{\boldsymbol \psi}_\mathfrak{F}(t)          \\[0.3em]
       \dot{\boldsymbol \psi}_{\mathfrak{S}}(t) 
     \end{bmatrix}=-\underbrace{\begin{bmatrix}
       \mathcal{L}_{g,d} & \mathcal{L}_{12}           \\[0.3em]
       \mathcal{L}_{21} & \mathcal{L}_{22}           
     \end{bmatrix}}_L\begin{bmatrix}
      {\boldsymbol \psi}_{\mathfrak{F}}(t)          \\[0.3em]
       {\boldsymbol \psi}_{\mathfrak{S}}(t) 
     \end{bmatrix}.
\label{eqn:mat}
\end{equation}
Given equation \eqref{eqn:fully}, we have that $ \mathcal{L}_{21}=0$ and $ \mathcal{L}_{22}=0 $. Hence the dynamics of the follower agents are given by 
\begin{align}
\dot{\boldsymbol \psi}_\mathfrak{F}(t) &= -\mathcal{L}_{g,d} \boldsymbol \psi_\mathfrak{F}(t) + \mathcal{L}_{12}\boldsymbol \psi_{\mathfrak{S}}(0).
\label{eqn:partial4}
\end{align}
 Here, $\mathcal{L}_{g,d}$ is the grounded Laplacian matrix, representing the interaction between the followers. When the graph is undirected, the grounded Laplacian is denoted by  $\mathcal{L}_{g,u}$. The submatrix $\mathcal{L}_{12}$ of the graph Laplacian captures the influence of the leaders on the followers. We make the following assumption in this paper. \newline
\textbf{Assumption 1:} In the directed graph $\mathcal{G}_d$, every follower can
be reached through a directed path from some leader. 

If Assumption 1 holds, the  states of the follower agents will converge to  some convex combination of the states of the leaders  \cite{Clarkbook}. Moreover, under Assumption 1, it is shown that the grounded Laplacian matrix is non-singular and $\lambda_1(\mathcal{L}_{g,d})$ is real and strictly positive and $\mathcal{L}_{g,d}^{-1}$ is a non-negative matrix \cite{cao}. In addition to the nominal dynamics \eqref{eqn:partial4}, we assume that there exists some disturbances (or perturbations) in the communications between the followers. In particular, consider the updating rule of each follower agent $v_j \in \mathfrak{F}$  is affected by a disturbance signal  $w_j(t)$ which turns \eqref{eqn:partial4} into  
\begin{align}
\dot{\boldsymbol {\psi}}_\mathfrak{F}(t) &= -\mathcal{L}_{g,d}{\boldsymbol \psi}_\mathfrak{F}(t)+ \mathcal{L}_{12}\boldsymbol \psi_{\mathfrak{S}}(0) + \boldsymbol w(t),\nonumber \\
\boldsymbol z(t)&={\psi}_\mathfrak{F}(t),
\label{eqn:padgrtial43}
\end{align}
where $\boldsymbol z(t)$ is the (full state) measurement.  Here $\boldsymbol w(t)$ is a vector representing the disturbances. We assumed that all followers are prone to be affected by the disturbances while the leaders are unaffected by the disturbances, since  they do not update their state. The objective is to quantify the effect of the external disturbance signals on the state of the follower agents. We use the system $\mathcal{H}_{\infty}$ norm of the transfer function $G(s)= (sI+\mathcal{L}_{g,d})^{-1}$ from $w(t)$ to $\boldsymbol z(t)$ defined as
\begin{align}
||G||_{\infty} &\triangleq \sup_{\omega\in \mathbb{R}}{\sigma_{n-|\mathfrak{S}|}\left(G(j\omega)\right)}.
\end{align}
Since $\mathcal{L}_{g,d}$ is Hurwitz, the above norm is finite.  We refer to $||G||_{\infty,u}$ and $||G||_{\infty,d}$ as system $\mathcal{H}_{\infty}$ norms when the underlying graph is undirected or directed, respectively. 
Before discussing the system norm of \eqref{eqn:padgrtial43}, we present the following definitions. 

\begin{definition}\textbf{(Positive Systems)}:
A linear system is called (internally) positive if and only if its state and output are non-negative for every non-negative input and every non-negative initial state.
\end{definition}
\begin{thm}[\cite{Rantzer}]
A continuous linear system $(A, B, C)$ is positive if and only if $A$ is a Metzler-matrix and $B$ and  $C$ are non-negative element-wise.  Moreover, for such a positive system with transfer function $G(s)=C(sI-A)^{-1}B$, the system $\mathcal{H}_{\infty}$ norm is obtained from the DC gain of the system, i.e.,  $||G||_{\infty}=\sigma_{n}(G(0))$, where  $\sigma_{n}$ is the maximum singular value of matrix $G(0)$.
\label{thm:metzlertheroem}
\end{thm}
It is clear that the evolution of follower agents \eqref{eqn:padgrtial43} together with full state measurements form a positive system. According to Theorem \ref{thm:metzlertheroem}, the system $\mathcal{H}_{\infty}$ norm from external disturbances to states of followers is $||G||_{\infty}=\sigma_{n-|\mathfrak{S}|}(\mathcal{L}_{g,d}^{-1})=\frac{1}{\sigma_{1}(\mathcal{L}_{g,d})}$. Hence, characterizing the system $\mathcal{H}_{\infty}$ norm  of \eqref{eqn:padgrtial43} is equivalent to determining the smallest singular value (or the eigenvalue for undirected networks) of the grounded Laplacian matrix, $\sigma_1(\mathcal{L}_{g,d})$.

\section{Bounds on System $\mathcal{H}_{\infty}$ Norms for directed  and undirected networks}\label{sec:boundssss}

In this section, we discuss upper and lower bounds for the smallest eigenvalue and singular value of the grounded Laplacian matrix. Then, based on these bounds, we provide graph-theoretic necessary and sufficient conditions for dynamics \eqref{eqn:padgrtial43} to have sufficiently small system norm.

\subsection{Undirected Network}

The following Theorem, which is an improved version of Theorem 1 in \cite{ArxiveRobutness}, provides some graph-theoretic bounds on $\lambda_1(\mathcal{L}_{g,u})$ for undirected networks. We use this result later in characterizing $\mathcal{H}_{\infty}$ performance in directed networks.

\begin{thm}\label{thm:dobare}
Consider a connected undirected graph $\mathcal{G}_u=\{\mathcal{V},\mathcal{E}\}$ with a set of leaders $\mathfrak{S} \subset \mathcal{V}$.  Let $\mathcal{L}_{g,u}$ be the grounded Laplacian matrix induced by $\mathfrak{S}$, and for each $v_i \in \mathfrak{F}$, let $\Gamma_i$ be the number of leaders in follower $v_i$'s neighborhood. Then
\begin{multline}
\max\left\{ \overbrace{\frac{1}{\mathcal{C}(v_k)}}^{(i)}, \overbrace{\Gamma_{\rm min}}^{(ii)}, \overbrace{\left(\frac{|\partial \mathfrak{S}|}{n-|\mathfrak{S}|}\right)x_{\rm min}}^{(iii)}\right\} \leq \lambda_1(\mathcal{L}_{g,u})
\leq  \min_{\emptyset \ne X\subseteq \mathcal{V\setminus \mathfrak{S}}} \frac{|\partial X|}{|X|} \leq \frac{|\partial \mathfrak{S}|}{n-|\mathfrak{S}|}\leq \Gamma_{\rm max},
\label{eqn:maineq}
\end{multline}
where $\mathcal{C}(v_k)$ is the closeness centrality of $v_k$ (sum of the shortest paths from all vertices to $v_k$) for any leader $v_k\in \mathfrak{S}$ and $x_{\rm min}$ is the smallest eigenvector component of $\mathbf{x}$, a non-negative eigenvector corresponding to $\lambda_1(\mathcal{L}_{g,u})$.\footnote{Eigenvector $\mathbf{x}$ is  normalized such that its largest component is $x_{\rm max}=1$.}
\end{thm}

\begin{figure}[t!]
\centering
\includegraphics[scale=.45]{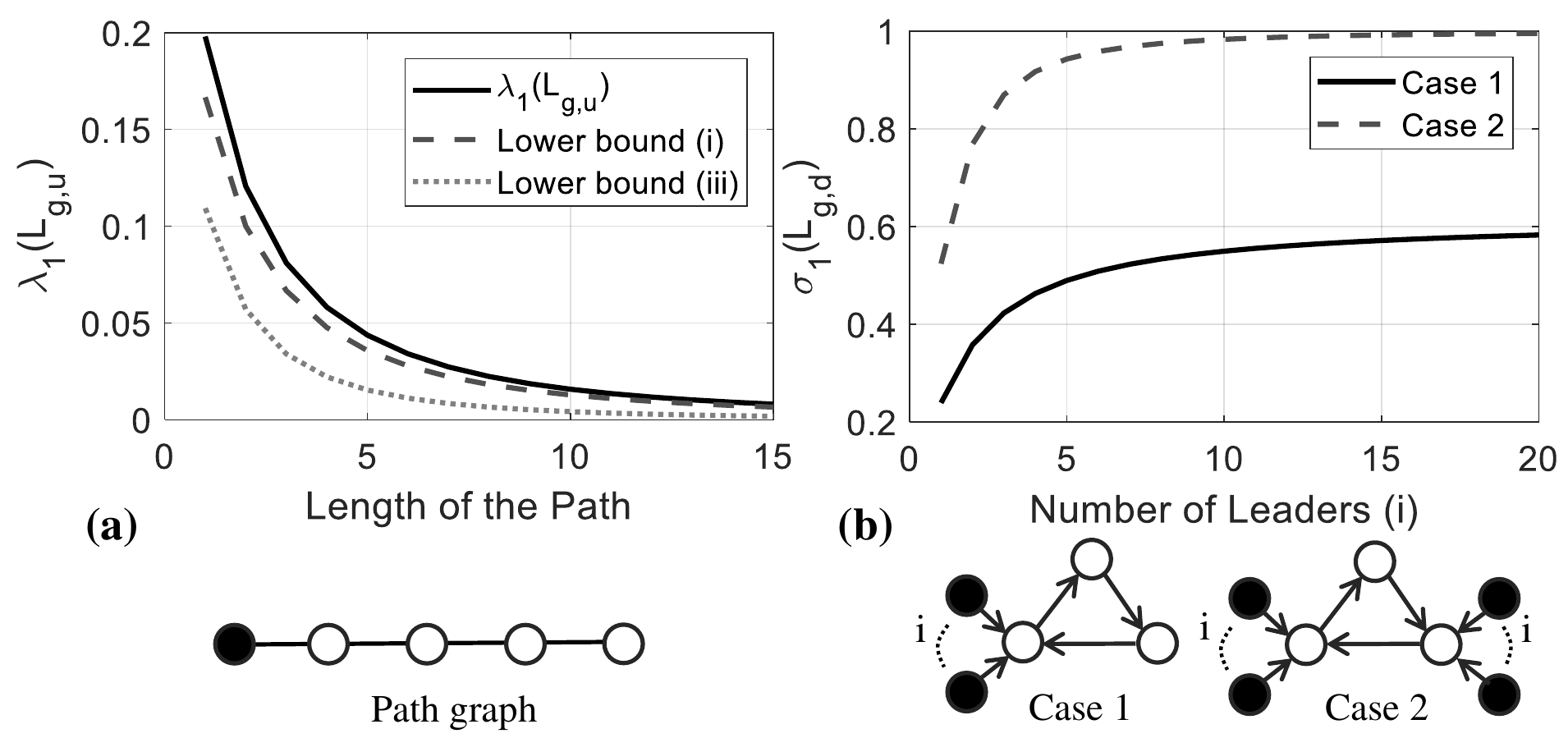}
\caption{(a) Tightness of lower bound $\frac{1}{\mathcal{C}(v_k)}$ in \eqref{eqn:maineq} compared to other lower bounds for a path graph, (b) Example showing the tightness of sufficient condition \eqref{eqn:suffcond}.}
\label{fig:etdsdvvnj}
\end{figure}

\begin{pf}
The proofs of all abounds exist in \cite{PiraniSundaramArxiv, ArxiveRobutness} except the lower bound $\frac{1}{\mathcal{C}(v_k)}$.
To show this lower bound, we pick any leader $v_k$ from the leader set and write the grounded Laplacian as $\mathcal{L}_{g,u}=\mathcal{L}_{g_k,u}+E$, where $\mathcal{L}_{g_k,u}$ is the grounded Laplacian for the graph where $v_k$ is the only leader and $E$ is a diagonal matrix showing the effect of the rest of the leaders. We know that $\mathcal{L}_{g_k,u}^{-1}$ is a non-negative matrix and due to the Perron-Frobenius theorem we have
\begin{equation}\label{eqn:perrr}
\lambda_{n-1}(\mathcal{L}_{g_k,u}^{-1})=\frac{1}{\lambda_1(\mathcal{L}_{g_k,u})}\leq \max_i[\mathcal{L}_{g_k,u}^{-1}]_i
\end{equation}
where $[\mathcal{L}_{g_k,u}^{-1}]_i$ is the $i$-th row sum of $\mathcal{L}_{g_k,u}^{-1}$. Based on the fact that 
\begin{align}
    [\mathcal{L}_{g_k,u}^{-1}]_{ii}&\leq {\rm dist}(v_i,v_k) \nonumber \\
    [\mathcal{L}_{g,u}^{-1}]_{ij}&\leq \min \{{\rm dist}(v_i,v_k), {\rm dist}(v_j,v_k)\},
\end{align}
where ${\rm dist}(v_i,v_j)$ is the length of the shortest path between $v_i$ and $v_j$ \cite{Miekkala}. By summing all elements of the $i$-th row of $\mathcal{L}_{g_k,u}^{-1}$ and substituting into \eqref{eqn:perrr} we get ${\lambda_1(\mathcal{L}_{g_k,u})}\geq \frac{1}{\mathcal{C}(v_k)}$. Now since $\mathcal{L}_{g,u}=\mathcal{L}_{g_k,u}+E$ and $E$ is a positive semidefinite matrix, based on Weyl's inequality, we conclude that by adding more leaders we have    $\lambda_1(\mathcal{L}_{g,u})\geq \lambda_1(\mathcal{L}_{g_k,u})$, and the result is obtained. \hfill\(\Box\)
\end{pf}
The newly proposed lower bound  $\frac{1}{\mathcal{C}(v_k)}$ on $\lambda_1(\mathcal{L}_g)$ is sometimes tighter than the others, as shown in Fig.~\ref{fig:etdsdvvnj} (a) for a path graph with a leader in one end. Here $\Gamma_{\rm min}=0$ which gives a trivial lower bound. We will revisit this bound in Section \ref{sec:platoonn}.

\subsection{Directed Network}

We first mention the following relation between the smallest eigenvalue of the grounded Laplacian and $\sigma_{1}(\mathcal{L}_{g,d})$. 
\begin{proposition}
For the network of leaders and followers satisfying Assumption 1, we have 
$\sigma_1(\mathcal{L}_{g,d}) \leq \lambda_1(\mathcal{L}_{g,d})$.
\label{prop:somebounds}
\end{proposition}
\begin{pf}
We know that for any square matrix $A_{n\times n}$, the spectral radius $\rho(A)$ is less than or equal to the largest singular value of $A$ \cite{Horn}, i.e., $\rho(A)\leq \sigma_{n}(A)$. Moreover, based on Assumption 1, we know that $\lambda_1(\mathcal{L}_{g,d})$ is real. Hence, we have 
$$
\frac{1}{\sigma_1(\mathcal{L}_{g,d})}=\sigma_{n-|\mathfrak{S}|}(\mathcal{L}_{g,d}^{-1})\geq \rho(\mathcal{L}_{g,d}^{-1})=\frac{1}{\lambda_1(\mathcal{L}_{g,d})},
$$
which yields the result. \hfill\(\Box\)
\end{pf}
Clearly for symmetric graphs, the smallest eigenvalue and singular value are the same and the bound of Proposition \ref{prop:somebounds} turns to equality. The problem becomes more challenging when the interaction network is non-symmetric. In this paper we investigate the relation between $\sigma_1(\mathcal{L}_{g,d})$ and $\lambda_1(\mathcal{L}_{g,u})$, and consequently the $\mathcal{H}_{\infty}$ norms in directed and undirected graphs, for specific classes of network. Before that, similar to what was discussed for undirected networks in Theorem \ref{thm:dobare}, we present the following theorem which provides some graph-theoretic bounds for the smallest singular value of the grounded Laplacian matrix.

\begin{thm}
Let $\mathcal{G}_d$ be a directed graph with leader set $\mathfrak{S}$ and suppose that Assumption 1 holds.  Then we have 
\begin{align}
\max\left\{0,\min_{1\leq i \leq n-|\mathfrak{S}|}\left\{\frac{1}{2}\left(\Delta_i-\delta_i+\Gamma_i \right) \right\}\right\}\leq \sigma_1(\mathcal{L}_{g,d}) \leq  \min \left\{\underbrace{\left(\frac{1}{n-|\mathfrak{S}|}\sum_{i=1}^{n-|\mathfrak{S}|}\Gamma_i^2\right)^{\frac{1}{2}}}_{(i)}, \underbrace{\left(\min_{v_i\in \mathcal{V}\setminus \mathfrak{S}}\{\Delta_i^2+\delta_i\}\right)^{\frac{1}{2}}}_{(ii)}
\right\}
\label{eqn:boundsonsep}
\end{align}
\label{thm:smallestsingularvalue}
\end{thm}
\begin{pf}
The lower bound comes from a Gershgorin-type bound proposed in \cite{Johnsoncharels}, in which for the smallest singular value of any $n$-by-$m$ matrix $A$ we have 
$$\min_{1\leq i \leq n}\left\{|a_{ii}|- \frac{1}{2}\left(\sum_{j=1,j\neq i}|a_{ij}|+\sum_{j=1,j\neq i}|a_{ji}| \right)\right\}\leq \sigma_1(A).$$
For the case where $A=\mathcal{L}_{g,d}$ we  have  $|a_{ii}|=\Delta_i$, $\sum_{j=1,j\neq i}|a_{ij}|=\Delta_i-\Gamma_i$, and $\sum_{j=1,j\neq i}|a_{ji}|=\delta_i$. For upper bounds, since $\mathcal{L}_{g,d}^T\mathcal{L}_{g,d}$ is a Hermitian matrix, from the Rayleigh quotient inequality \cite{Horn}, we have
$$
\sigma_1 \leq \left(\boldsymbol y^T \mathcal{L}_{g,d}^T\mathcal{L}_{g,d} \boldsymbol y \right)^{\frac{1}{2}},
$$
for all $\boldsymbol y\in \mathbb{R}^{n-|S|}$ with $\boldsymbol y^T\boldsymbol y=1$. The upper bound  $ \left(\frac{1}{n-|\mathfrak{S}|}\sum_{i=1}^{n-|\mathfrak{S}|}\Gamma_i^2\right)^{\frac{1}{2}}$ is then obtained by choosing $\boldsymbol y=\frac{1}{\sqrt{n-|\mathfrak{S}|}}\mathbf{1}_{n-|\mathfrak{S}|}$. The upper bound $\left(\min_{v_i\in \mathcal{V}\setminus \mathfrak{S}}\{\Delta_i^2+\delta_i\}\right)^{\frac{1}{2}}$ is obtained by choosing $\boldsymbol y=\mathbf{e}_i$ where $i=\arg \min_{v_i\in \mathcal{V}\setminus \mathfrak{S}} \{\Delta_i^2+\delta_i\}$.  \hfill\(\Box\)
\end{pf}

Note that none of the upper bounds provided in Theorem \ref{thm:smallestsingularvalue} are in general greater or less than the other one. In order to show the tightness of those bounds, we have the following example.
\begin{exmp}
In Fig.~\ref{fig:etfdvnj} (a),  both (non-zero) lower and upper bounds are equal to 1 and \eqref{eqn:boundsonsep} is tight. Computing the lower bound for  graph (b) gives $\frac{1}{2}.$ The upper bounds (i) and (ii) give $\frac{1}{\sqrt{2}}\approx 0.7$ and 1, respectively, while its actual value is $ \sigma_1(\mathcal{L}_{g,d})=0.62$. In graph (c), the lower bound gives 0.5 and upper bounds (i) and (ii) give $\frac{5}{2}$ and $2$, respectively. It should be noted that in graph (b),  upper bound (i) is tighter, while in graph (c), upper bound (ii) is tighter. 
\begin{figure}[t!]
\centering
\includegraphics[scale=.5]{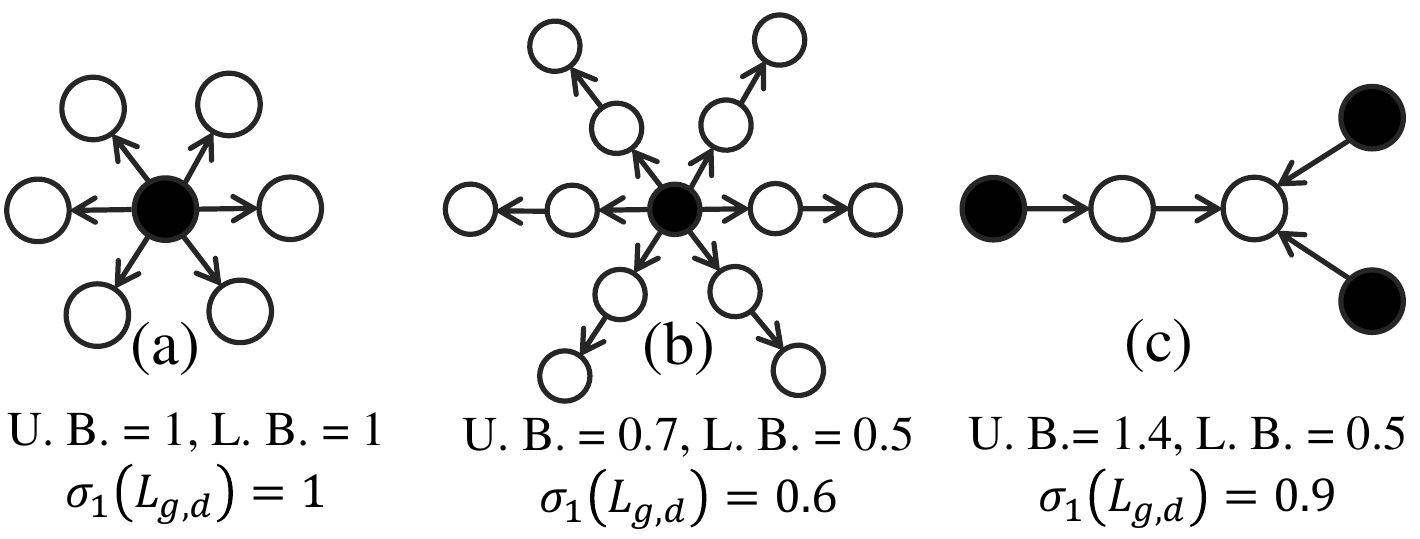}
\caption{Upper and lower bounds in \eqref{eqn:boundsonsep} and the value of $\sigma_1(\mathcal{L}_{g,d})$  for different networks.}
\label{fig:etfdvnj}
\end{figure}
\end{exmp}

 Theorem \ref{thm:smallestsingularvalue} yields the following corollary which presents graph-theoretic conditions for dynamics \eqref{eqn:padgrtial43} to have sufficiently small $\mathcal{H}_{\infty}$ norm.
\begin{corollary}
A sufficient condition for dynamics \eqref{eqn:padgrtial43} to have
$||G||_{\infty}\leq \gamma$ is to have
\begin{equation}\label{eqn:suffcond}
    \min_{v_i\in \mathcal{V}\setminus \mathfrak{S}}\{\Delta_i-\delta_i+\Gamma_i\}\geq \lceil\frac{2}{\gamma}\rceil.
\end{equation}
Moreover, for dynamics \eqref{eqn:padgrtial43} to have $||G||_{\infty}\leq \gamma$ it is necessary to satisfy both conditions
\begin{align}
\sum_{i=1}^{n-|\mathfrak{S}|}\Gamma_i^2 \geq \lfloor\frac{n-|\mathfrak{S}|}{\gamma^2}\rfloor, \quad
\min_{v_i\in \mathcal{V}\setminus \mathfrak{S}}\{\Delta_i^2+\delta_i\}\geq \lfloor\frac{1}{\gamma^2}\rfloor,
\end{align}
and based on the fact that $\sum_{i=1}^{n-|\mathfrak{S}|}\Gamma_i^2\leq (n-|\mathfrak{S}|)\Gamma_{\rm max}^2$, another (simpler) necessary condition for $||G||_{\infty}\leq \gamma$ is to have $\Gamma_{\rm max} \geq \lfloor\frac{1}{\gamma}\rfloor$.
\label{cor:chemidd}
\end{corollary}

The conditions mentioned in Corollary \ref{cor:chemidd} can give some clues in designing networks with desired robustness. 
For instance, based on sufficient condition \eqref{eqn:suffcond} we should have $ \min_{v_i\in \mathcal{V}\setminus \mathfrak{S}}\{\Delta_i-\delta_i+\Gamma_i\}\geq 2$ to ensure that the system is non-expansive, i.e., $||G||_{\infty}\leq 1$. Fig.~\ref{fig:etdsdvvnj} (b) shows  a directed cycle of size 3 which does not satisfy sufficient condition \eqref{eqn:suffcond} for $||G||_{\infty}\leq 1$. Here,  even if we add many leaders to one or two of the nodes (Cases 1 and 2 in Fig.~\ref{fig:etdsdvvnj} (b)) we can not get $\sigma_1(\mathcal{L}_{g,d})\geq 1$, i.e., $||G||_{\infty}\leq 1$. However, for sufficient condition \eqref{eqn:suffcond} to satisfy we need to add two leaders to each node to get $||G||_{\infty}\leq 1$. Moreover, based on the necessary condition $\Gamma_{\rm max} \geq \lfloor\frac{1}{\gamma}\rfloor$, we know that it is impossible to get $||G||_{\infty}<1$ with only one leader, regardless of the interconnections within the follower nodes. Moreover, based on $\min_{v_i\in \mathcal{V}\setminus \mathfrak{S}}\{\Delta_i^2+\delta_i\}\geq \lfloor\frac{1}{\gamma^2}\rfloor$, it is impossible to get $||G||_{\infty}<1$ when there exists a leaf in the graph (a node with $\Delta_i=1$ and $\delta_i=0)$. 

The following proposition shows the effect of the number of leaders on the system robustness. 

\begin{proposition}\label{prop:addingleader}
Consider a directed graph $\mathcal{G}_d$ with leader set $\mathfrak{S}$  which satisfies Assumption 1. Then adding leaders to the network does not increase the system $\mathcal{H}_{\infty}$ norm. 
\end{proposition}
\begin{pf}
We show that by adding extra leaders, $\sigma_1(\mathcal{L}_{g,d})$ does not decrease. By adding leaders we can write $\mathcal{L}_{g,d}=\hat{\mathcal{L}}_{g,d}+E$, where $\hat{\mathcal{L}}_{g,d}$ is the grounded Laplacian of the original graph with $|\mathfrak{S}|$ leaders and $E$ is the diagonal matrix which shows the effect of the extra leaders. Then we have
\begin{align}
\lambda_1(\mathcal{L}_{g,d}^T\mathcal{L}_{g,d})=\lambda_1(\hat{\mathcal{L}}_{g,d}^T\hat{\mathcal{L}}_{g,d}+\hat{\mathcal{L}}_{g,d}^TE+E^T\hat{\mathcal{L}}_{g,d}+E^TE)
\geq \lambda_1(\hat{\mathcal{L}}_{g,d}^T\hat{\mathcal{L}}_{g,d}).
\label{eqn:extraleadersss}
\end{align}
 The inequality in \eqref{eqn:extraleadersss} is due to the Weyl's inequality and comes from the fact that matrices $E^TE$ and $\hat{\mathcal{L}}_{g,d}^TE+E^T\hat{\mathcal{L}}_{g,d}$ are positive semidefinite (the latter has a weighted Laplacian structure). \hfill\(\Box\)
\end{pf}
The effect of adding extra leaders on increasing $\sigma_1(\mathcal{L}_{g,d})$ is shown in Fig.~\ref{fig:etdsdvvnj} (b) as well. In the following section, we discuss the relation between system $\mathcal{H}_{\infty}$ norms in directed networks and their undirected counterparts. Moreover, we discuss some inconsistencies between directed and undirected networks in the sense of the behaviour of the system $\mathcal{H}_{\infty}$ norm in response of adding extra edges.

\section{Relations Between $\mathcal{H}_{\infty}$ norm of directed and undirected networks}\label{sec:relations}

In this section, we compare the system $\mathcal{H}_{\infty}$ norms in directed graphs and their undirected counterparts. Since there is no specific relation between the two cases in general graphs, we focus on particular classes of networks for which we can derive explicit expressions for the relation between system $\mathcal{H}_{\infty}$ norms in directed and undirected networks.  Hence, we restrict the attention to two classes of networks, namely {\it balanced digraphs} and {\it directed trees}.

\subsection{Balanced Digraphs}
\label{sec:influence}
Balanced Digraphs are directed graphs for which the in-degree  and out-degree of each node are equal. For these graphs, we will show that the system $\mathcal{H}_{\infty}$ norm of a directed network is no worst than twice of that of undirected network. Before that, we have the following lemma. 

\begin{lem}[\cite{dingara2}]
The system $\mathcal{H}_{\infty}$ norm of a positive system with asymmetric interactions is upper bounded by the $\mathcal{H}_{\infty}$ norm of the symmetric parts of the dynamic matrix. 
\label{lem:abb}
\end{lem}

\begin{thm}
Consider a directed graph $\mathcal{G}_d$ with leader set $\mathfrak{S}$  which satisfies Assumption 1. If the subgraph of the follower agents is balanced, then the system $\mathcal{H}_{\infty}$ norm of \eqref{eqn:padgrtial43} satisfies $\|G\|_{{\infty, d}} \leq 2||G||_{\infty , u}$.
\label{thm:brl12}
\end{thm}
\begin{pf}
According to Lemma \ref{lem:abb}, we have
\begin{equation}
    \|G\|_{{\infty, d}} \leq \frac{2}{|\lambda_{1}(\mathcal{L}_{g,d}+\mathcal{L}_{g,d}^{\sf T})|}\,.
    \label{eqn:brl}
\end{equation}
Moreover, for balanced graphs  we have $\mathcal{L}_{g,d}+\mathcal{L}_{g,d}^T=\bar{\mathcal{L}}_g+E$ for some diagonal and positive semidefinite matrix $E$, where $\bar{\mathcal{L}}_g$ is the grounded Laplacian matrix corresponding to the undirected network. Thus based on Weyl's inequality we have  $\lambda_{1}(\bar{\mathcal{L}}_g)\leq \lambda_{1}(\mathcal{L}_{g,d}+\mathcal{L}_{g,d}^T)$, and together with \eqref{eqn:brl} the result is obtained.  \hfill\(\Box\)
\end{pf}

The following example shows that one can not modify Theorem \ref{thm:brl12} to get  $\|G\|_{{\infty, d}} \leq ||G||_{\infty , u}$ in balanced digraphs.

\begin{figure}[t!]
\centering
\includegraphics[scale=.45]{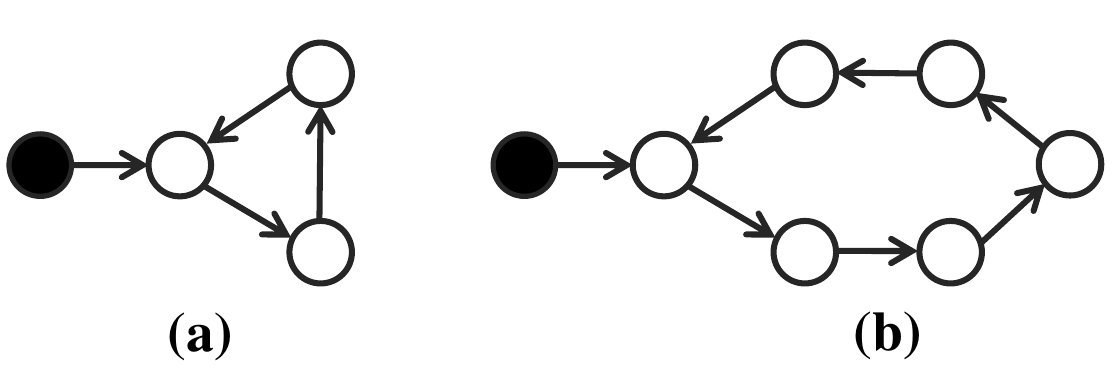}
\caption{Balanced Digraphs: (a) example showing $||G||_{\infty , d}>||G||_{\infty , u}$, (b) example  showing $||G||_{\infty , d}<||G||_{\infty , u}$.}
\label{fig:etnj}
\end{figure}

\begin{exmp}
 As shown in Fig.~\ref{fig:etnj} (a), which is a balanced graph of followers, we have  $||G||_{\infty , u}=3.73$ and $||G||_{\infty , d}=4.18$. Moreover, if we increase the length of the loop from 3 to 6, Fig.~\ref{fig:etnj} (b), we have $||G||_{\infty , u}=9.19>8.85=||G||_{\infty , d}$. This shows that the bound proposed in Theorem \ref{thm:brl12} is tight and $||G||_{\infty , d}<||G||_{\infty , u}$ does not always hold.
\end{exmp}

We should note that, as shown in \cite{Dingharajovanovic}, changing the direction of the edges in a balanced digraph does not change the system $\mathcal{H}_{\infty}$ norm.

\subsection{Directed Trees}
In this section, we focus on directed networks whose undirected counterparts are  trees, i.e., connected graphs without cycles. In the following theorem, we discuss the relation between the system $\mathcal{H}_{\infty}$ norm of \eqref{eqn:padgrtial43}  in directed and undirected trees. 
\begin{thm}
Consider a directed graph $\mathcal{G}_d$ with leader set $\mathfrak{S}$  which satisfies Assumption 1. If the subgraph of the followers is a tree, the system $\mathcal{H}_{\infty}$ norm of \eqref{eqn:padgrtial43} satisfies
\begin{equation}
    \frac{1}{\min_{i\in \mathcal{V}\setminus \mathfrak{S}} \Delta_i}\leq ||G||_{\infty , d}\leq ||G||_{\infty , u}^{\frac{1}{2}},
\end{equation}
with $||G||_{\infty , d}=||G||_{\infty , u}^{\frac{1}{2}}$ if there exists a single leader, i.e., $|\mathfrak{S}|=1$. 
\label{thm:treeing}
\end{thm}
\begin{pf}
For the case where there exists a single leader in the network, for each follower node there is exactly one incoming edge, as otherwise a cycle will be made in the underlying undirected graph. If we write the grounded Laplacian matrices of directed and undirected graphs with a single leader by $\hat{\mathcal{L}}_{g,d}$ and $\hat{\mathcal{L}}_{g,u}$, respectively, we know that
$\hat{\mathcal{L}}_{g,d}$
is triangular with diagonal elements 1. In this case we have
$\hat{\mathcal{L}}_{g,d}^T\hat{\mathcal{L}}_{g,d}= \hat{\mathcal{L}}_{g,u}$. It is due to the fact that each diagonal element of $\hat{\mathcal{L}}_{g,d}^T\hat{\mathcal{L}}_{g,d}$ is $\Delta_i^2+\delta_i=1+\delta_i$ which is equal to the degree of each node in the undirected network. Moreover, the off-diagonal elements of $\hat{\mathcal{L}}_{g,d}^T\hat{\mathcal{L}}_{g,d}$ are just the reflections of the lower-triangle elements of $\hat{\mathcal{L}}_{g,d}$ to the upper-triangle. 
Hence, we have $\lambda_1(\hat{\mathcal{L}}_{g,d}^T\hat{\mathcal{L}}_{g,d})=\lambda_1(\hat{\mathcal{L}}_{g,u})$ which shows that the $\mathcal{H}_{\infty}$ norm for the undirected graph is equal to the square of its directed version, i.e. $||G||_{\infty , u}=||G||_{\infty , d} ^2$. For $|\mathfrak{S}|>1$, similar to Proposition \ref{prop:addingleader} we can write the grounded Laplacians as ${\mathcal{L}}_{g,d}=\hat{\mathcal{L}}_{g,d}+E$ and ${\mathcal{L}}_{g,u}=\hat{\mathcal{L}}_{g,u}+E$, where $E$ shows the effect of the rest of the leaders. Considering the fact that $E^TE \succeq E$, i.e., $E^TE-E$ is positive semidefnite, and $\hat{\mathcal{L}}_{g,d}^TE+E^T\hat{\mathcal{L}}_{g,d}$ is also positive semidefnite, we get
\begin{align}
\lambda_1(\mathcal{L}_{g,d}^T\mathcal{L}_{g,d})=\lambda_1(\hat{\mathcal{L}}_{g,d}^T\hat{\mathcal{L}}_{g,d}+\hat{\mathcal{L}}_{g,d}^TE+E^T\hat{\mathcal{L}}_{g,d}+E^TE)
\geq \lambda_1(\hat{\mathcal{L}}_{g,u}+E)=\lambda_1({\mathcal{L}}_{g,u}),
\label{eqn:extraleadersss}
\end{align}
which yields $||G||_{\infty, d}^2\leq ||G||_{\infty, u}$.

For the lower bound, via an appropriate permutation of rows, matrix $\mathcal{L}_{g,d}$ can be put into a triangular form. Then we have $\lambda_1(\mathcal{L}_{g,d})=\min_{i\in \mathcal{V}\setminus \mathfrak{S}} \Delta_i$ (the minimum in-degree in the subgraph of followers) and according to proposition \ref{prop:somebounds}, we have $||G||_{\infty , d}\geq \frac{1}{\min_{i\in \mathcal{V}\setminus \mathfrak{S}} \Delta_i}$. 
\hfill\(\Box\)
\end{pf}

\begin{exmp}
 As shown in Fig.~\ref{fig:etvnj}, the system $\mathcal{H}_{\infty}$ norm of the directed tree with multiple leaders can be larger or smaller than that of the undirected graph, depends on if $||G||_{\infty , d}\geq 1$ or not. In this figure, for graph (a) we have $||G||_{\infty,d}=1.14<1.7=||G||_{\infty,u}$ and for graph (b) we have $||G||_{\infty,d}=0.54>0.5=||G||_{\infty,u}$. The lower bound  mentioned in Theorem \ref{thm:treeing} is $ \frac{1}{\min_{i\in \mathcal{V}\setminus \mathfrak{S}} \Delta_i}=1$ for graph (a) and $ \frac{1}{\min_{i\in \mathcal{V}\setminus \mathfrak{S}} \Delta_i}=0.5$ for graph (b) which shows that it is close to the actual value of $||G||_{\infty , d}$. 

\begin{figure}[t!]
\centering
\includegraphics[scale=.45]{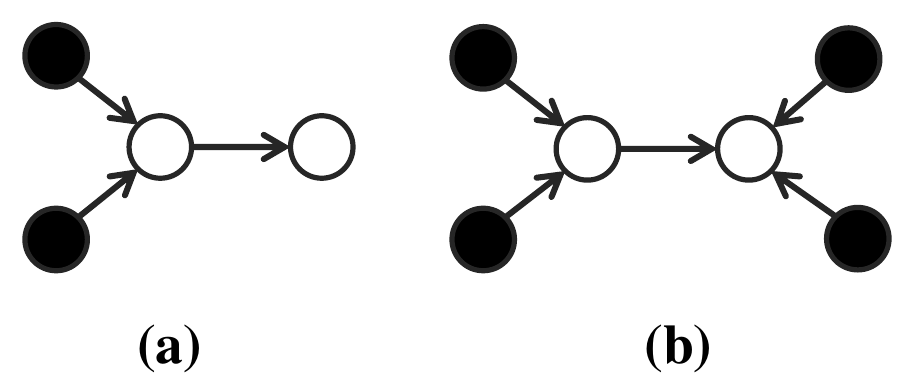}
\caption{Trees: (a) example showing $||G||_{\infty,d}<||G||_{\infty,u}$, (b) example showing $||G||_{\infty,d}>||G||_{\infty,u}$.}
\label{fig:etvnj}
\end{figure}
\end{exmp}

Combining Theorem \ref{thm:dobare} with Theorem \ref{thm:treeing} yields the following corollary.
\begin{cor}
Consider a directed graph $\mathcal{G}_d$ with leader set $\mathfrak{S}$  which satisfies Assumption 1. If the subgraph of the followers is a tree, the system $\mathcal{H}_{\infty}$ norm of \eqref{eqn:padgrtial43} satisfies
\begin{equation}\label{eqn:upperclose}
||G||_{\infty, d}\leq \min_{v\in \mathfrak{S}}\mathcal{C}(v)^{\frac{1}{2}},
\end{equation}
where $\mathcal{C}(v)$ is the closeness centrality of any leader node $v$, as mentioned in Theorem \ref{thm:dobare}.
\end{cor}

We will come back to the above result in Section \ref{sec:platoonn}.

\subsection{Effect of Adding Edges to Directed Trees}

\begin{figure}[t!]
\centering
\includegraphics[scale=.60]{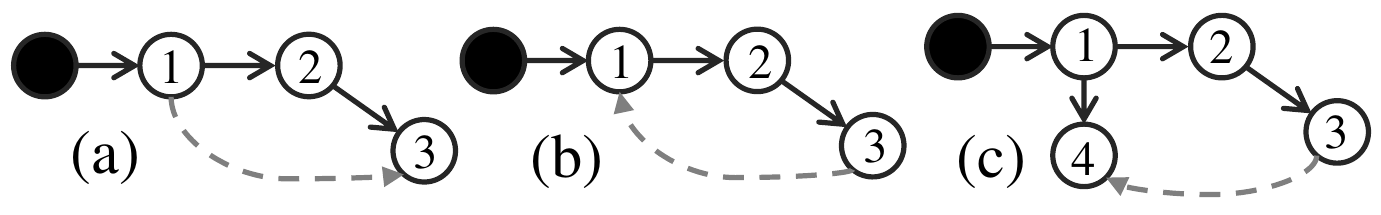}
\caption{(a), (b) Adding an edge to a directed path with opposite directions, (c) Adding an edge to a directed tree.}
\label{fig:etvsdfnj}
\end{figure}

In the previous subsection, we discussed directed graphs whose undirected counterpart is a tree. In this subsection, we consider the effect of adding extra directed edges to a directed tree on the system $\mathcal{H}_{\infty}$ norm. We present an observation in Fig.~\ref{fig:etvsdfnj}.  
For graph (a), an additional directed edge (grey dashed line) is added to a directed path which does not make a directed cycle; however, for case (b) the additional edge makes a cycle.  The system $\mathcal{H}_{\infty}$ norm of \eqref{eqn:padgrtial43} for the path graph (before adding the directed edge) is $||G||_{\infty,d}=2.25$, for case (a) is $||G||_{\infty,d}=1.99$ and for case (b) is $||G||_{\infty,d}=4.18$. Hence, it implies that adding a directed edge to a directed path to make a cycle deteriorates the $\mathcal{H}_{\infty}$ performance and adding an edge to a path which does not make a cycle improves the performance. This observation is intuitive, since adding a cycle to the network results in  the information (and uncertainties) to circulate (and thus propagate) in a part of the network. However, the opposite is not true for general trees, as shown in graph (c). In particular, for general trees, the $\mathcal{H}_{\infty}$ performance can be deteriorated even for an edge addition which does not make a cycle in the network. For graph (c) in  Fig.~\ref{fig:etvsdfnj}, before adding a directed edge we have  $||G||_{\infty,d}=2.40$ and after adding that (which does not make a cycle)  we have $||G||_{\infty,d}=2.56$. Based on these observations, we present a result in Theorem \ref{thm:directedtreerobustness}. Before that, we need the following definition and Lemma \ref{lem:vanmeigh}.

\begin{definition}\textbf{(Leader-Rooted path and Interfering  Edges):}
A {\it leader-rooted-path} in a directed tree is a directed path which starts from a leader's neighbor and ends at one of the leaf nodes, i.e., nodes with $\Delta_i=1$ and $\delta_i=0$. 
Moreover, given a leader-rooted-path labeled by $v_1, v_2, ..., v_{\ell}$, where $v_1$ is the leader's neighbor and $v_{\ell}$ is the leaf, the additional two directed edges  to this path, named $(v_i,v_j)$ and $(v_k,v_h)$, are called {\it interfering} if $j>k$. 
\end{definition}

\begin{figure}[t!]
\centering
\includegraphics[scale=.58]{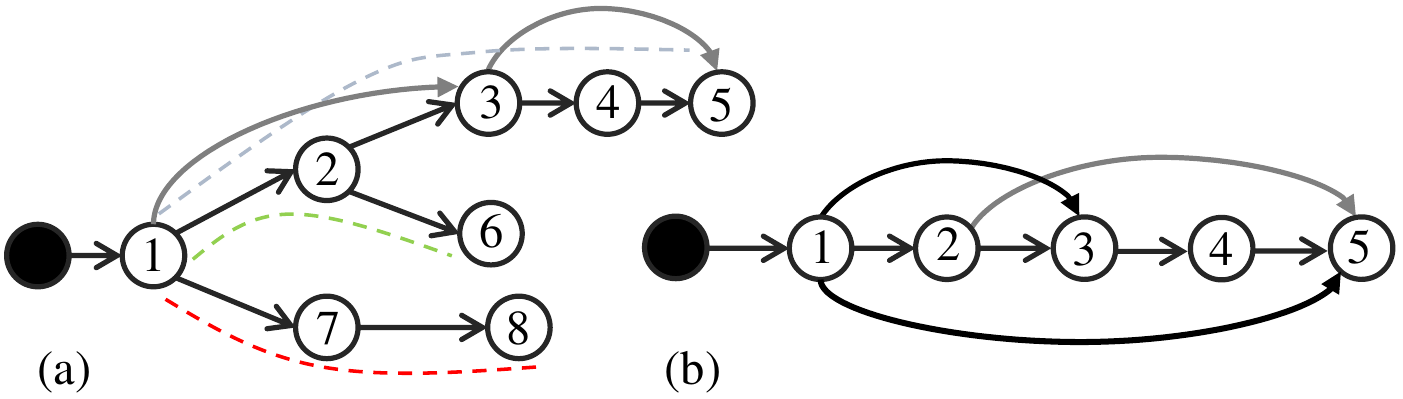}
\caption{(a) Example of a tree with three leader-rooted paths, (b) Example of additional interfering paths. }
\label{fig:leaderrooted}
\end{figure}

Three different leader-rooted-paths in a directed tree are shown in Fig.~\ref{fig:leaderrooted} (a) by dashed lines. Moreover, the two additional edges (in grey) are not interfering. However, in Fig.~\ref{fig:leaderrooted} (b), the grey and the black edges are interfering. 
Figure \ref{fig:leaderrooted} (b) is an example which shows that even for a path graph, if the additional edges are interfering, then they may increase the $\mathcal{H}_{\infty}$ norm. In this example, before adding the grey edge, we have $||G||_{\infty}=2.65$ and after that we have $||G||_{\infty}=2.66$. Based on the above observations and definitions, we will present Theorem \ref{thm:directedtreerobustness}. Before that, we have the following lemma. 
\begin{lem}[\cite{vanmeighem}]\label{lem:vanmeigh}
If one element in a non-negative matrix $A$ is increased, then the largest
eigenvalue is also increased. The increase is strict for irreducible matrices.
\end{lem}

According to Lemma \ref{lem:vanmeigh}, we can conclude that if one element in a non-negative matrix $A$ is decreased (but  still be positive), then the largest eigenvalue is also decreased. Based on this, we present the following theorem.
\begin{thm}\label{thm:directedtreerobustness}
Consider a directed graph $\mathcal{G}_d$ with leader set $\mathfrak{S}$  which satisfies Assumption 1. If the subgraph of the followers is a tree, then adding non interfering directed edges between the nodes in a leader-rooted-path which does not make a cycle will decrease the system $\mathcal{H}_{\infty}$ norm and adding  a directed edge between two nodes in a leader-rooted-path which makes a cycle will increase the system $\mathcal{H}_{\infty}$ norm. 
\end{thm}
\begin{pf}
When an edge is added from node $j$ to node $i$, we can write the new Laplacian matrix as $\Tilde{\mathcal{L}}_{g,d}= \mathcal{L}_{g,d} + \mathbf{e}_i\mathbf{e}^T_{ij}$, where $\mathbf{e}_i$ is a vector which is $1$ in $i -th$ place and zero elsewhere and $\mathbf{e}_{ij}=\mathbf{e}_i-\mathbf{e}_j$. We know that $\Tilde{\mathcal{L}}_{g,d}^{-1}$ is a non-negative matrix. Without loss of generality, we label the nodes in the leader-rooted-path containing $i$ and $j$ as $1, 2, ..., j, ..., i, ..., r$, where $1$ belongs to the leader's neighbor and $r$ is the length of the path.  First, we will show that adding a directed edge from $j$ to $i$ decreases some elements of $\Tilde{\mathcal{L}}_{g,d}^{-1}$ and adding an edge from $i$ to $j$ increases them. If we use Sherman-Morrison formula \cite{shermanmorrison} we get
\begin{align}\label{eqn:vnosup}
\Tilde{\mathcal{L}}_{g,d}^{-1}=\left( \mathcal{L}_{g,d} + \mathbf{e}_i\mathbf{e}^T_{ij} \right)^{-1}=\mathcal{L}_{g,d}^{-1}-\frac{\mathcal{L}_{g,d}^{-1}\mathbf{e}_i\mathbf{e}^T_{ij}\mathcal{L}_{g,d}^{-1}}{1+\mathbf{e}^T_{ij}\mathcal{L}_{g,d}^{-1}\mathbf{e}_i}.
\end{align}
As both nodes $i$ and $j$ are in the same leader-rooted-path, the block of matrix ${\mathcal{L}}_{g,d}^{-1}$ from row $1$ to row $r$ and column $1$ to column $r$ is in the form of a lower triangular matrix whose lower triangle elements are all 1 (it can be easy verified by solving the corresponding block in ${\mathcal{L}}_{g,d}^{-1}$ from ${\mathcal{L}}_{g,d}{\mathcal{L}}_{g,d}^{-1}=I$).  Hence, $\mathcal{L}_{g,d}^{-1}\mathbf{e}_i\mathbf{e}^T_{ij}\mathcal{L}_{g,d}^{-1}$ is a non-negative matrix and $\mathbf{e}^T_{ij}\mathcal{L}_{g,d}^{-1}\mathbf{e}_i=1$.  Thus, based on \eqref{eqn:vnosup}, it implies that  $\Tilde{\mathcal{L}}_{g,d}^{-1}$ is non-negative and its elements are not larger than those of ${\mathcal{L}}_{g,d}^{-1}$. Likewise, the elements of $\Tilde{\mathcal{L}}_{g,d}^{-T}\Tilde{\mathcal{L}}_{g,d}^{-1}$ are not larger than the elements of ${\mathcal{L}}_{g,d}^{-T}{\mathcal{L}}_{g,d}^{-1}$ and based on Lemma \ref{lem:vanmeigh} we conclude that $\lambda_{n-|\mathfrak{S}|}({\mathcal{L}}_{g,d}^{-T}{\mathcal{L}}_{g,d}^{-1}) \geq  \lambda_{n-|\mathfrak{S}|}(\Tilde{\mathcal{L}}_{g,d}^{-T}\Tilde{\mathcal{L}}_{g,d}^{-1})$ which proves the claim. For the case where the additional edge makes a cycle, i.e., from $i$ to $j$,  with the similar argument we can show that the elements of $\frac{\mathcal{L}_{g,d}^{-1}\mathbf{e}_i\mathbf{e}^T_{ij}\mathcal{L}_{g,d}^{-1}}{1+\mathbf{e}^T_{ij}\mathcal{L}_{g,d}^{-1}\mathbf{e}_i}$ are non positive and the rest of the proof is similar.

If another edge, called $(v_x, v_y)$, is added to the same leader-rooted path which is not interfering with edge $(v_j, v_i)$, since ${\mathcal{L}_{g,d}^{-1}\mathbf{e}_i\mathbf{e}^T_{ij}\mathcal{L}_{g,d}^{-1}}$ affects columns $j+1$ to $i$ and we have $x\geq i$, then  the block of matrix ${\mathcal{L}}_{g,d}^{-1}$ from row $x$ to row $y$ and column $x$ to column $y$ is not affected by the previous edge $(v_j, v_i)$ and the rest of the proof is the same.
\hfill\(\Box\)
\end{pf}

 Theorem \ref{thm:directedtreerobustness} proved a general scenario which includes observations in Fig.~\ref{fig:etvsdfnj} (a) and (b). Clearly, in graph shown in Fig.~\ref{fig:etvsdfnj} (c) nodes 3 and 4 are not in the same leader-rooted path; hence, the sufficient condition mentioned in  Theorem \ref{thm:directedtreerobustness} can not be applied to that case. 
 
 \section{Application: Vehicle platooning }\label{sec:platoonn}
 
 In this section, we apply the results from previous sections to a vehicle platooning problem. 
 Consider a connected network of $n$ vehicles. The position and longitudinal velocity of each vehicle $v_i$ is denoted by scalars $p_i$ and $u_i$, respectively. Each vehicle $v_i$ is able to communicate with its neighbor vehicles and transfer its kinematic parameters, e.g., velocity. \footnote{Transmitting vehicle's states such as velocity is common in standard short-range vehicular communications \cite{IEEE1609}.} 
 The objective for each follower vehicle  is to track a reference velocity  $u^*$. This  desired reference velocity is calculated  by a reference vehicle  in order to optimize the fuel consumption \cite{van2015fuel}. The dynamics of vehicle $v_i$ is governed by $\ddot{p}_i(t)= q_i(t)$, or in vector notation
\begin{equation}
\ddot{\mathbf{p}}(t)=\mathbf{q}(t)+\mathbf{w}(t),
\label{eqn:law}
\end{equation}
where  $\mathbf{p}(t)=[{p}_1(t), {p}_2(t), ..., {p}_n(t)]^T$ is the vector of positions and  $\mathbf{q}(t)$ is the vector of the control law and $\mathbf{w}(t)$ represents the  disturbance effect through inter-vehicular communications. The following control laws for each follower and reference vehicle are considered, \cite{hao2013stability}
\begin{align}
{q}_i(t)=
  \begin{cases}
    \sum_{j\in \aleph_i}(u_j(t)-u_i(t))       & \quad  \forall  v_i\in \mathfrak{F},\\
   0  & \quad \forall v_i \in \mathfrak{S}.\\
  \end{cases}
\label{eqn:single}
\end{align}

\begin{figure}[t!]
\centering
\includegraphics[scale=.45]{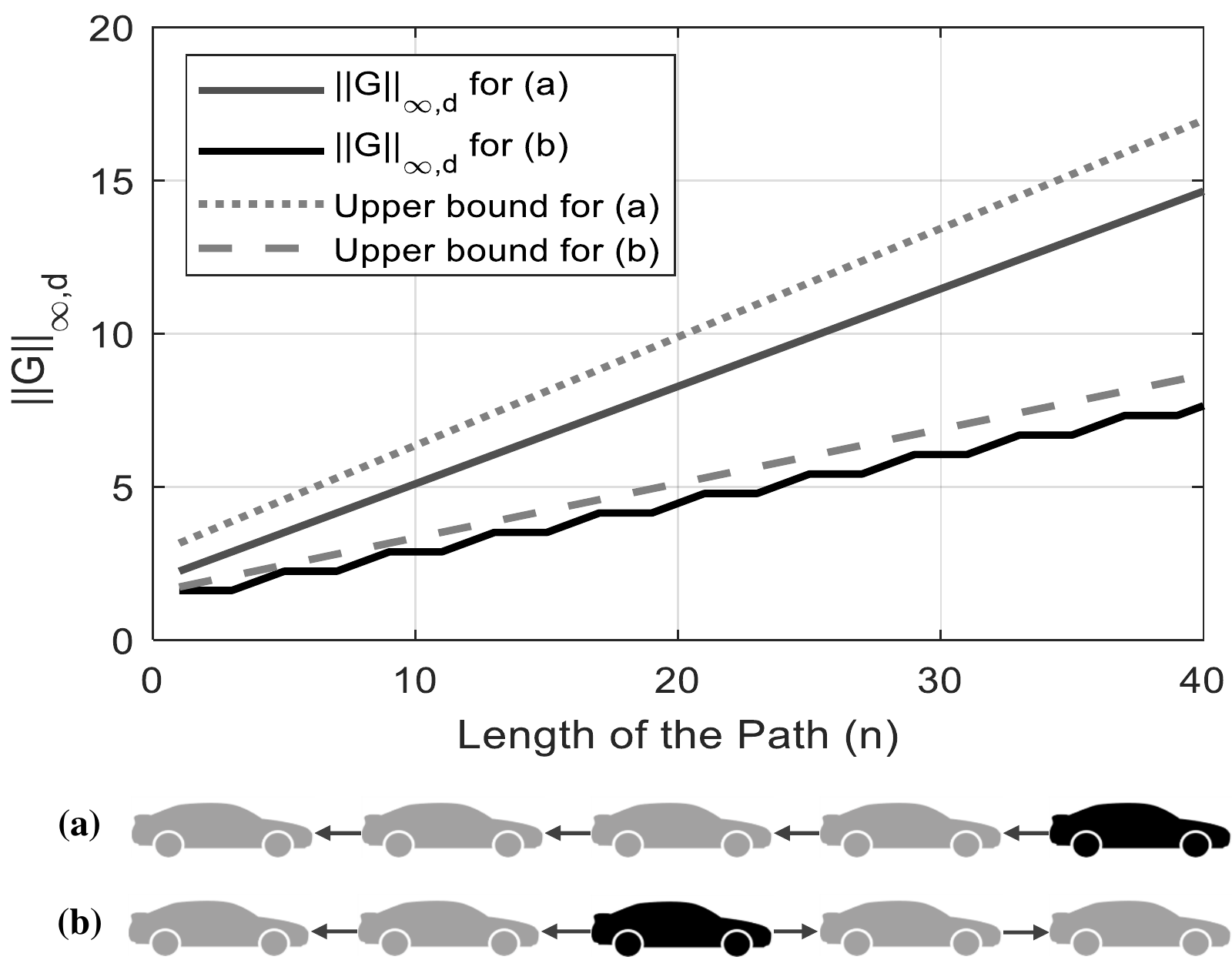}
\caption{Vehicle Platoons, showing  the effect of the leader's location on the $\mathcal{H}_{\infty}$ performance of the system.}
\label{fig:etdsdvvsdnj}
\end{figure}

It is common to use system $\mathcal{H}_{\infty}$ norm to quantify the effect of communication noises/disturbances in vehicle platooning \cite{seiler2004disturbance, hao2013stability}. Here, we discuss the effect of the leader's position on the $\mathcal{H}_{\infty}$ performance of reference velocity tracking via applying the results of the previous sections. As shown in Fig.~\ref{fig:etdsdvvsdnj}, for a platoon of connected cars on a directed path graph, we set the leader either in one end of the platoon, Fig.~\ref{fig:etdsdvvsdnj} (a), or in the middle of the platoon, Fig.~\ref{fig:etdsdvvsdnj} (b). For both cases, the system $\mathcal{H}_{\infty}$ norm of \eqref{eqn:law} from $\mathbf{w}(t)$ to $\dot{\mathbf{p}}(t)$, the vector of velocities, for various platoon lengths are shown in Fig.~\ref{fig:etdsdvvsdnj}. Moreover, the upper bounds predicted by \eqref{eqn:upperclose}, which was a function of the closeness centrality of the leader, are shown in this figure for each scenario. According to this result, placing the reference vehicle in an appropriate position in the network has a considerable impact on the robustness of the platoon to communication disturbances.   

\begin{figure}[t!]
\centering
\includegraphics[scale=.415]{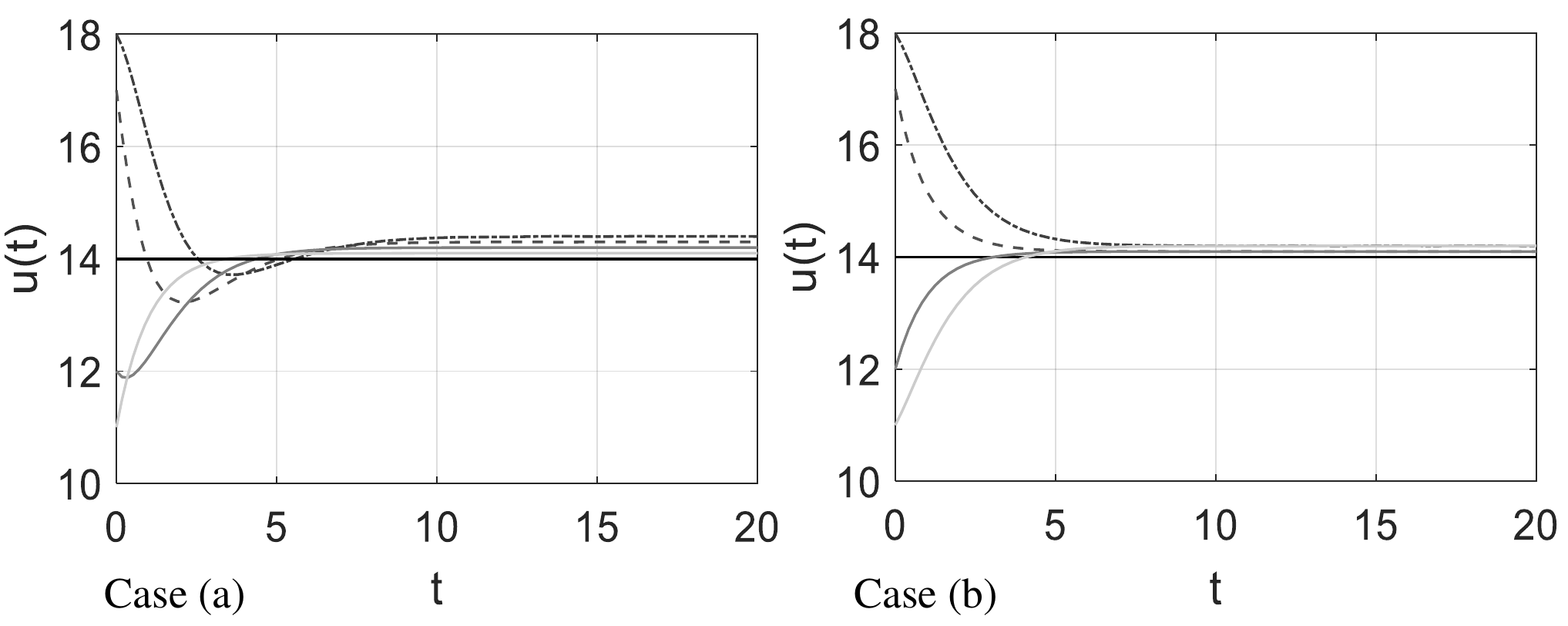}
\caption{Time-responses of velocity tracking scenario for cases (a) and (b) in Fig.~\ref{fig:etdsdvvsdnj}. }
\label{fig:etdsdvvsfgdnj}
\end{figure}

This effect can be seen in the time-response of the system as well. Figure \ref{fig:etdsdvvsfgdnj} shows the time response of the velocity tracking scenario for cases (a) and (b) in Fig.~\ref{fig:etdsdvvsdnj} (the reference velocity is $u^*=14  \frac{m}{s}$), where a constant additive disturbance is added to the dynamics, i.e., $\mathbf{w}(t)=0.1\times\mathbf{1}$ in \eqref{eqn:law}. According to this figure, the steady-state error for case (b) is less than that of case (a) and since the $L_2$-norm of the disturbance signal for two cases are equal, it results in smaller $L_2$ gain of the system (which is equal to $\mathcal{H}_{\infty}$ norm for LTI systems \cite{Boyyd}). 

\section{Summary}
 \label{sec:conclusion}
In this paper, a graph-theoretic approach to the $\mathcal{H}_{\infty}$ performance of leader following consensus dynamics on directed graphs was studied. The relation between the system $\mathcal{H}_{\infty}$ norm for directed and undirected networks for specific classes of graphs, i.e., balanced digraphs and directed trees, was discussed. Moreover, the effects of adding directed edges to a directed tree on the resulting system $\mathcal{H}_{\infty}$ norm was investigated. At the end, the results were applied to  a reference velocity tracking problem in a platoon of  vehicles. A future avenue for further research is to generalize the comparison made between $\mathcal{H}_{\infty}$ norm for directed and undirected networks to more diverse topological structures. Moreover, extending these results to the $\mathcal{H}_{\infty}$ performance of second-order consensus dynamics on directed networks, which has direct applications in vehicle formation problems, is an interesting direction of research.

\appendix

\bibliographystyle{plain}        
\bibliography{refs}         



\end{document}